\documentclass[10pt]{article}  
\usepackage{amsmath}
\usepackage{amssymb}
\usepackage{dsfont}
\usepackage{epsfig}
\usepackage{graphics}
\usepackage{theorem}

\theoremheaderfont{\scshape}

\theorembodyfont{\normalfont}

\numberwithin{equation}{section}

\title{Models for cohesive sediments describing the evolution of the characteristics of particles}
\author{
Emmanuel Fr\'enod \thanks{Universit\'e Europ\'enne de Bretagne, Lab-STICC (UMR CNRS 3192), 
Universit\'e de Bretagne-Sud, Centre Yves Coppens, Campus de Tohannic,
F-56017, Vannes} }

\newcommand{\nit}{\mathbb{N}}
\newcommand{\rit}{\mathbb{R}}

\newcommand{\pos}{\mathbf{x}}
\newcommand{\vit}{\mathbf{v}}
\newcommand{\vvit}{\mathbf{V}}
\newcommand{\etur}{\mathbf{k}}
\newcommand{\tdr}{\mathbf{\varepsilon}}
\newcommand{\ph}{\textit{p}\textrm{H}}

\newcommand{\rhc}{\rho}
\newcommand{\rhca}{\rho^a}
\newcommand{\rhcb}{\rho^b}
\newcommand{\rhct}{\tilde\rho}
\newcommand{\rhch}{\hat\rho}

\newcommand{\ssvel}{{{\cal W}_{\hspace{-1pt}s}}}

\newcommand{\mass}{{m}}
\newcommand{\ffl}{{\cal F}}
\newcommand{\Ba}{{{\cal B}_{\hspace{-1pt}a}}}
\newcommand{\Bf}{{{\cal B}_{\hspace{-1.5pt} f}}}
\newcommand{\Be}{{{\cal B}_{\hspace{-1pt}e}}}
\newcommand{\Bet}{{\tilde{\cal B}_{\hspace{-1pt}e}}}
\newcommand{\De}{{{\, \cal D}_{\hspace{-2pt}eq}}}
\newcommand{\tr}{{{\, \cal T}_{\hspace{-3pt}eq}}}
\newcommand{\G}{{{\bf G}}}
\newcommand{\Gov}{{\overline{\bf G}}}
\newcommand{\f}{{{f}}}
\newcommand{\Ld}{{{L}}}
\newcommand{\Lr}{{{L}^{-1}}}
\newcommand{\Bfl}{{\overline{\bf B}^{~}_{\hspace{-1pt}\textit{lf}}}}
\newcommand{\Bfli}{{{\bf B}^{\hspace{2pt}i}_{\hspace{-1pt}\textit{lf}}}}
\newcommand{\Bla}{{\overline{\bf B}^{~}_{\hspace{-1pt}\textit{la}}}}
\newcommand{\Blai}{{{\bf B}^{\hspace{2pt}ij}_{\hspace{-1pt}\textit{la}}}}
\newcommand{\Bgf}{{\overline{\bf B}^{~}_{\hspace{-1.5pt}\textit{gf}}}}
\newcommand{\Bgfi}{{{\bf B}^{\hspace{2pt}ij}_{\hspace{-1.5pt}\textit{gf}}}}
\newcommand{\Bga}{{\overline{\bf B}^{~}_{\hspace{-1.5pt}\textit{ga}}}}
\newcommand{\Bgai}{{{\bf B}^{\hspace{2pt}ikl}_{\hspace{-1.5pt}\textit{ga}}}}

\newcommand{\lm}{{\lambda_\textrm{min}}}
\newcommand{\lb}{{\lambda_\textrm{bio}}}
\newcommand{\mm}{{{\cal M}_\textrm{min}}}
\newcommand{\mb}{{{\cal M}_\textrm{bio}}}
\newcommand{\nd}{{{\cal N}_{\hspace{-1pt}d}}}

\newcommand{\rrr}{{r}}
\newcommand{\roum}{{p}}

\newcommand{\ldt}{{\tilde\lambda}}
\newcommand{\ldtp}{{\tilde\lambda}'\hspace{0.5pt}}

\newcommand{\lnpp}{{(\lambda^d-\lambda'^d)^{1/d}}}
\newcommand{\lnss}{{(\lambda^d-\lambda''^d)^{1/d}}}
\newcommand{\lpnp}{{(\lambda'^d-\lambda^d)^{1/d}}}
\newcommand{\lptnpt}{{(\ldtp^d-\ldt^d)^{1/d}}}
\newcommand{\lnptpt}{{(\ldt^d-\ldtp^d)^{1/d}}}
\newcommand{\lpnpd}{{\lambda'^d-\lambda^d}}
\newcommand{\lnppm}{{(\lambda^d-\lambda'^d)^{(1-d)/d}}}
\newcommand{\lnssm}{{(\lambda^d-\lambda''^d)^{(1-d)/d}}}
\newcommand{\lnssminv}{{(\lambda^d-\lambda''^d)^{(d-1)/d}}}
\newcommand{\lpnpm}{{(\lambda'^d-\lambda^d)^{(1-d)/d}}}
\newcommand{\lptnptm}{{(\ldtp^d-\ldt^d)^{(1-d)/d}}}
\newcommand{\lnptptm}{{(\ldt^d-\ldtp^d)^{(1-d)/d}}}

\newcommand{\lnppsum}{{(\lambda^d+\lambda'^d)^{1/d}}}

\newcommand{\ds}{\displaystyle}

\newcommand{\fracp}[2]{\frac{\partial #1}{\partial #2}}

\begin{document}
\maketitle

{\small {\bf Abstract - } 
The goal of this paper is to set up a framework designed to take into account the characteristics of sediment particles when transported by water. 
Our protocol consists in describing the characteristics of sediment particles via an additional 
variable, and to build operators involving this new variable, modeling the evolution of the particle characteristics.
Several such operators are proposed, some based on principles of relaxation toward an equilibrium, 
and others on a description of the particles' aggregation and fragmentation process.
A discrete version of the latter is also offered for numerical settings.
}

~

{\small {\bf Keywords - } 
Modeling, 
Cohesive Sediments, 
Relaxation Models,
Aggregation and Fragmentation Models,
Integro-Differential Equations.
}

{\small
\tableofcontents}

\section{Introduction}
\label{intr}
In view of the evolution of the climate, and the increasingly stringent requirements 
in terms of feasibility and impact studies before dredging or building sea walls
or harbors, the estuary morpho-dynamical issue is becoming a topic of major
significance.
As a consequence, the behavior of the complex sediments found in estuaries 
is today widely studied, analyzed and modeled. 
The research effort concerns every aspect, from measurement protocols
to simulations of deposition, erosion, transport by water, wave action,
turbulence results and flocculation processes. 

In all these aspects,  modeling has a key role to play. 
To summarize, the modeling of cohesive sediments involves three compartments, as well as 
the interactions which link them together. 
The first compartment deals with fluid field forecasts, and involves Navier-Stokes 
or Shallow Water type equations, possibly involving 
turbulence and using a propagating eddy viscosity. 
The second describes sediment behavior when it is deposited on the seabed, and the 
third models the transport of sediment particles when suspended in the water
column. 
Roughly speaking, this last model propels sediment particles at the same velocity as the 
fluid added to a settling velocity. 
It may also take into account the action of turbulence, using the eddy viscosity as a dissipative 
effect on sediment particles.

~

The present paper is situated in this context, and focuses on the issue of the transport 
of sediment particles by the water column. 
It offers a robust and flexible framework which takes into account the 
evolution and alterations of the sediment particle characteristics (size, mass, porosity,
 \emph{etc.}) while the particles are in the water column, concomitantly with other 
 phenomena (transport, settling, turbulence). 
The main idea consists in introducing a mass density of sediment particles
$\rhc(t,\pos,\lambda)$ depending on time $t$, position $\pos =(x,y,z)$, and also an
additional variable $\lambda \in \Lambda$, which describes the particles' characteristics.

~

To present this idea, this article will begin by giving a proper definition of
$\rhc$, in the context of this framework, consisting in a two-part integro-differential
equation. The first part consists in a differential operator acting on 
$\rhc$ and describing the action of water as it transports particle. The second part 
is an integral operator, named $\G$, modeling the evolution of the particles' characteristics.   
Examples of sets $\Lambda$  of characteristics,  variables $\lambda$ and 
integral operators $\G$ will then be given, before exploring the way in which existing 
aggregation models may be translated into the present framework.
Finally, a discrete instantiation of set $\Lambda$ and operator $\G$ are given for numerics.

\section{Guiding ideas}
\label{GI}
The following elements will be taken for granted: 
A given estuary may be represented by a regular subset $\Omega \in \rit^3$, provided with 
coordinates $x,y,z$, where the $x-$axis is horizontal and points toward the east, the  $y-$axis is
horizontal and points toward the north, and the $z-$axis is vertical and points toward the sky.
At any time $t\in \rit^+$ and in any point  $\pos =(x,y,z)$ of $\Omega$, the water velocity 
$\vit$ with coordinates $(u,v,w)$ in $(x,y,z)-$coordinate system may be computed using
a Navier-Stokes-type system, possibly  involving eddy viscosity.
The salinity $S$ and the temperature $T$ may be obtained by solving 
advection-diffusion equations, possibly involving eddy viscosity.
Moreover, the energy of turbulence $\etur$ and its dissipation rate $\tdr$ may be 
computed using, for instance, a $\etur-\tdr$ model. The eddy viscosity involved in
the equation describing the evolution of $\vit$, $S$ and $T$ may be computed from
$\etur$ and $\tdr$. We may also suppose that the water $\ph$ and the amount of organic matter 
per liter of water $O$ are available. 
In the sequel, $\ffl$ denotes the fluid field
\begin{gather}
\label{GI1}
\ffl = (\vit, S,T, \etur,\tdr,\ph, O),
\end{gather}
which ranges in $\rit^3\times(\rit^+)^6$ and which depends on $t$ and $\pos$.

~

The main idea to be explored hereafter consists in assuming that the characteristics
of the sediment particles may be described by a variable $\lambda$ belonging to a given continuous
space $\Lambda$ and that, at time $t\in \rit^+$, the mass distribution of suspended matter of 
type $\lambda\in \Lambda$ and in point $\pos\in\Omega$ may be described by a measure, which is absolutely continuous with respect to the Lebesgue measure, and with density
$\rhc(t,\pos,\lambda)$.
The precise definition of $\rhc$ states that for any subset $\omega\subset \Omega\times\Lambda$,
the mass of sediment particles with position and characteristics situated within $\omega$ is
\begin{gather}
\label{GI2}  
\int_\omega \rhc(t,\pos, \lambda) \, d\pos d\lambda,
\end{gather}
at any given time $t$.

Actually, particles result from the assembly of elementary sediment particles. 
Hence, the characteristic of a given particle naturally belongs to a discrete space
$\tilde \Lambda$, and the mass distribution $R$, with respect to the characteristic
variable $\lambda$, is naturally a sum of Dirac mass
\begin{gather}
R = \sum_{\tilde \lambda \in \tilde \Lambda} 
N(t, \pos, \tilde \lambda) \; d\pos \; \delta_{\lambda =\tilde \lambda} .
\end{gather} 
Hence, when making the above assumption, we consider that, on the observation
scale, $R$ may be replaced (or approached) by $ \rhc(t,\pos, \lambda) \, d\pos d\lambda$.

From $\rhc$, the mass density of suspended matter at $t$ is defined by 
\begin{gather}
\label{GI3}  
\rrr(t,\pos) = 
\int_\Lambda \rhc(t,\pos, \lambda) \, d\lambda.
\end{gather}

The evolution of $\rhc$, over time, is supposed to be the result of transport by water, settling, 
diffusion by turbulence and aggregation, fragmentation and, more generally, the shape and mass 
evolution of the particles.
In point of fact, $\rhc$ should be seen as the solution to:
\begin{multline}
\label{GI4}  
\fracp{\rhc}{t} 
+ U(\ffl,\lambda)\fracp{\rhc}{x}
+ V(\ffl,\lambda)\fracp{\rhc}{y}
+ \big(W(\ffl,\lambda)-\ssvel(\lambda,\rrr)\big)\fracp{\rhc}{z}
\\
-\bigg(
\fracp{\Big(\mu(\ffl,\lambda)\ds\fracp{\rhc}{x}\Big)}{x}
+ \fracp{\Big(\mu(\ffl,\lambda)\ds\fracp{\rhc}{y}\Big)}{y}
+ \fracp{\Big(\nu(\ffl,\lambda)\ds\fracp{\rhc}{z}\Big)}{z}
\bigg)
= \G (\ffl,\rhc,\lambda).
\end{multline}
The first four terms in (\ref{GI4}) are the time derivative of $\rhc$ following the trajectories
induced by velocity $(U,V,W-\ssvel)$.
Velocity $\vvit = \vvit(\ffl,\lambda)=(U(\ffl,\lambda),V(\ffl,\lambda),W(\ffl,\lambda))$ is the velocity
transmitted by water to the particles. It depends on the fluid field $\ffl$ and on the particles' characteristics
$\lambda$. But in most situations, it is reasonable to set $ \vvit(\ffl,\lambda)= \vit$, meaning that 
water transmits its velocity directly to sediment particles. Settling velocity $\ssvel(\lambda,\rrr)$ is
the velocity at which particles fall toward the seabed.
It is natural to consider that $\ssvel$ depends on particle characteristics, with the idea that 
the heavier a particle, the faster it falls. Moreover, if the particle density $\rrr$ in the water 
column is high, settling may be slowed or hindered by the proximity of many particles.
In this case, $\ssvel$ has to depend on $\rrr$.
The fifth term of the equation's left-hand side conveys the fact that sediment particles undergo diffusion.
This involves a horizontal diffusion coefficient $\mu$ and a vertical diffusion coefficient $\nu$.
The diffusion phenomenon essentially comes from turbulence. Hence, choosing for 
$\mu(\ffl,\lambda)$ and $\nu(\ffl,\lambda)$ the water's usual eddy viscosity, given
by $\tilde c_\mu \frac{\etur^2}{\tdr}$ (where $\tilde c_\mu \sim 90$) is not completely 
unreasonable, at least for particle characteristics $\lambda$ corresponding to small sizes.
For particles characteristics corresponding to sizes bigger than the turbulent structures,
other choices must be made. 
The right-hand side of equation(\ref{GI4}) models the evolution of the particle
characteristics.
Clearly, the way in which particles combine, fragment, grow or, more generally evolve is linked
with the fluid field, especially aspects such as temperature, salinity, $\ph$, concentration of 
organic matter  and certainly turbulence energy.
Thus $\G$ depends on $\ffl$. The next section will provide examples of operator $\G$ and a
discussion of its properties.
\section{Relaxation models}
\label{RM}
This section gives simple examples of operator $\G$.

In the following paragraphs, the operator $\G$ is not built up from physical considerations, but only
by considering the asymptotic evolution of sediment particles when
environmental conditions, or in other words, when the fluid field $\ffl$, remains the same
over a long period of time.
\subsection{A mass-preserving relaxation model with one-dimensional $\Lambda$}
In this first example, $\Lambda$ is assumed to be $\rit^+$ and $\lambda \in\Lambda$
stands for particle size. Considering this certainly supposes that all the sediment 
particles under consideration have the same shape, and that they can be characterized
by a one-dimensional parameter. If the particles are one-dimensional, $\lambda$ is the
particle length, and particle mass is in direct proportion to $\lambda$.  
If they are two-dimensional, $\lambda$ is the particle diameter and particle mass
is in direct proportion to $\lambda^2$. If they are tridimensional, $\lambda$ is also
the diameter but particle mass is in proportion to $\lambda^3$.
\subsubsection{Operator building}
The operator-building process presented here is influenced by Bhatnagar, Gross \& Krook \cite{BGK},
who offered a kinetic model for gas dynamics.

If it is well established that, when fluid field $\ffl$ remains the same in a given
place over a long period, the mass density distribution with respect to $\lambda$ is given by the following equilibrium distribution function:
\begin{gather}
\label{RM1} 
\roum \De (\ffl, \lambda),
\end{gather}
where $\roum$ is a mass density with respect to $\pos-$variable ($\roum$ depends on  $\pos$)
and where $\De (\ffl, .)$ is a density probability defined by $\Lambda$  (in particular, it 
satisfies:
\begin{gather}
\label{RM2} 
\int_{\Lambda} \De (\ffl, \lambda) \, d \lambda =1, \text{ for all } \ffl \; )
\end{gather}
then, introducing a relaxation time $\tr$, $\G$ may be defined as
\begin{gather}
\label{RM3} 
\G(\ffl,\rhc,\lambda) = -\frac{1}{\tr}
\bigg( \rhc - \Big(\int_{\Lambda} \rhc \, d\lambda' \Big)\De(\ffl,\lambda)\bigg),
\end{gather}
or
\begin{gather}
\label{RM4} 
\G(\ffl(t,\pos),\rhc(t,\pos,.),\lambda) = -\frac{1}{\tr}
\bigg( \rhc(t,\pos,\lambda) -\Big( \int_{\Lambda} \rhc(t,\pos,\lambda') \, d\lambda' \Big)
\De(\ffl(t,\pos),\lambda)\bigg).
\end{gather}
\subsubsection{Example of function $\De$}
\begin{figure}
\centering  
\begin{tabular}{rrr}
\includegraphics[width=4cm]{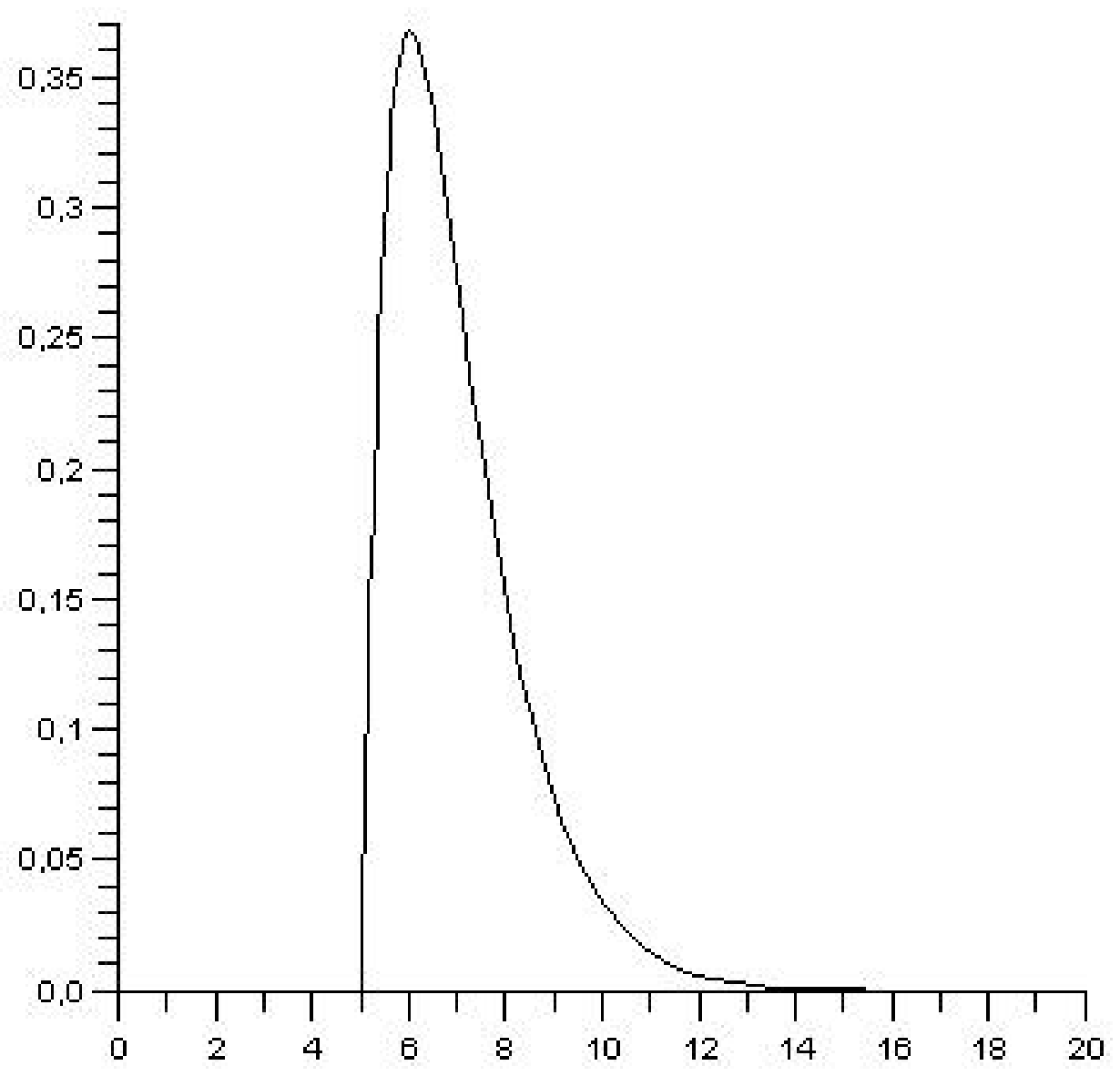}  & 
\includegraphics[width=4cm]{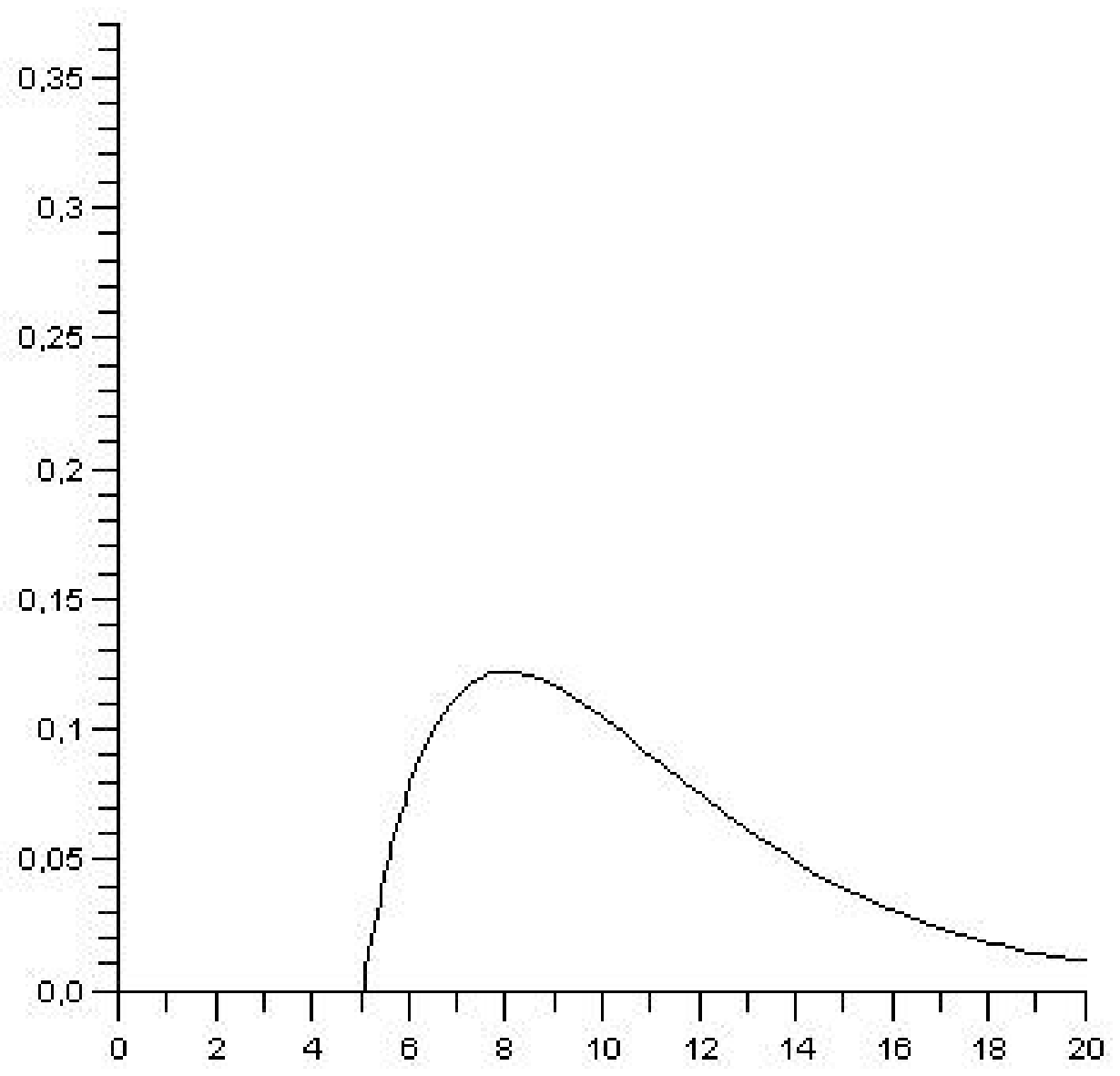}\vspace{-17pt}   
\\
 \small $ _\lambda$ ~~~~&   $ _\lambda$~~~~~\vspace{-100pt} 
 \\
\hspace{-35 pt} $ _{\De(\ffl,\lambda)}$ \hfill ~  & 
 \hspace{-35 pt} $ _{\De(\ffl,\lambda)}$ \hfill ~\hspace{40 pt}\vspace{75pt} 
\end{tabular}
\caption{Function $\De(\ffl,\lambda)$  defined by  (\ref{RM5}) for $\lm=5$, $\sigma(\ffl)=1$ (left)  and $\sigma(\ffl)=3$ (right)}
\label{Fig1} 
\end{figure}

As function $\De(\ffl,\lambda)$, defined for any $\ffl \in \rit^3\times (\rit^+)^6$ and any
$\lambda\in\rit^+$, we can choose:
\begin{gather}
\label{RM5} 
\begin{aligned}
&\De(\ffl,\lambda) & & = ~ 0 & &\text{ if }  \lambda< \lm \\
&                             & & 
        = ~ \frac{\lambda-\lm}{(\sigma(\ffl))^2}\:\exp \bigg( \frac{\lambda-\lm}{(\sigma(\ffl))}\bigg)
                                            & &\text{ otherwise },
\end{aligned}
\end{gather}
which is drawn in figure \ref{Fig1} for $\lm=5$ and for  $\sigma(\ffl)=1$ on the left
and $\sigma(\ffl)=3$ on the right. This choice makes it possible to take into account that particles
cannot be smaller than $\lm$, and that the variability of particle size, for a given
equilibrium, depends on a function $\sigma(\ffl)$ of the fluid field.

\subsubsection{Properties}
When used in (\ref{GI4}), the operator $\G$ defined by (\ref{RM3})  pushes
$\rhc(t,\pos,\lambda)$ toward $\big( \int_{\Lambda}\rhc(t,\pos,\lambda') \, d\lambda' \big)$ $\De(\ffl,\lambda)$ 
at any time and place, with a relaxation time of $\tr$. Moreover, the action of $\G$
in (\ref{GI4}) does not influence the evolution of the total mass of sediment. 

In order to be more precise, a function $\rhc(t,\lambda)$, not depending on $\pos$, 
which is solution to
\begin{gather}
\label{RM6}  
\fracp{\rhc}{t} = \G(\ffl,\rhc,\lambda),
\end{gather} 
for a fixed vector $\ffl \in\rit^3\times(\rit^+)^3$, has the following properties.
First, the quantity
\begin{gather}
\label{RM7}  
\int_{\Lambda} \rhc(t,\lambda) \, d\lambda,
\end{gather}
remains constant over time and, secondly, for every $\lambda\in\Lambda$, the quantity
\begin{gather}
\label{RM8}  
\rhc(t,\lambda) - \Big( \int_{\Lambda} \rhc(t,\lambda') \, d\lambda' \Big) \De(\ffl,\lambda),
\end{gather}
is divided by $e$ after any period of time of length $\tr$.

~

Property (\ref{RM7}) may be seen by integrating (\ref{RM6}):
\begin{multline}
\label{RM9}  
\fracp{\Big(\ds \int_{\Lambda} \rhc \, d\lambda\Big)}{t} 
= \int_{\Lambda}\fracp{\rhc}{t}\, d\lambda
= -\frac{1}{\tr} \int_{\Lambda}
\bigg( \rhc - \Big(\int_{\Lambda} \rhc(., \lambda' )\, d\lambda' \Big)\De(\ffl,\lambda)\bigg)\, d\lambda
\\
=  - \frac{1}{\tr} \bigg(\int_{\Lambda}  \rhc(., \lambda')\, d\lambda 
                                  - \int_{\Lambda} \rhc(., \lambda' ) \, d\lambda'  \bigg) =0.
\end{multline}
Property (\ref{RM8}) may be seen by computing the solution to equation (\ref{RM6}) leading, 
for any $t>s$, to
\begin{gather}
\label{RM10}  
\rhc(t, \lambda) - \Big(\int_{\Lambda}  \rhc(t, \lambda')\, d\lambda'\Big) \De(\ffl,\lambda)
=\bigg( \rhc(s, \lambda) - \Big( \int_{\Lambda}  \rhc(s, \lambda')\, d\lambda'\Big) \De(\ffl,\lambda) \bigg)
e^{(s-t)/\tr}.
\end{gather}
\subsection{ A non-mass preserving relaxation model with one-dimensional $\Lambda$}
\begin{figure}
\centering 
\begin{tabular}{rr}
\includegraphics[width=3.5cm]{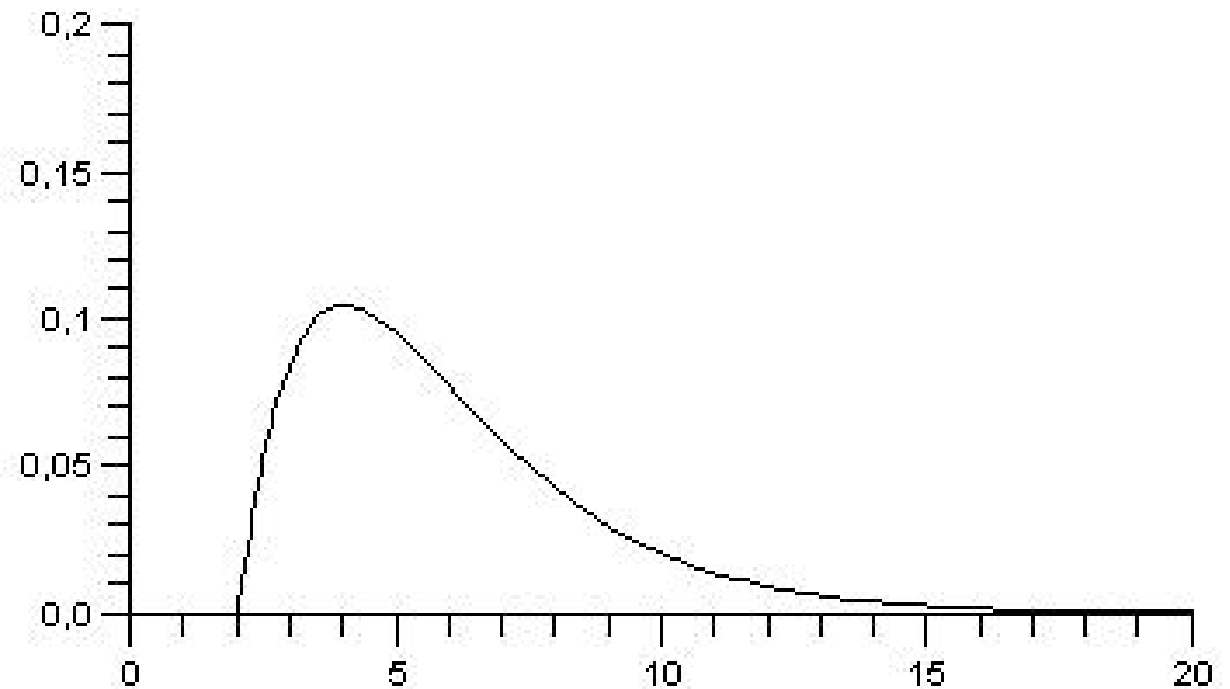} &
\includegraphics[width=3.5cm]{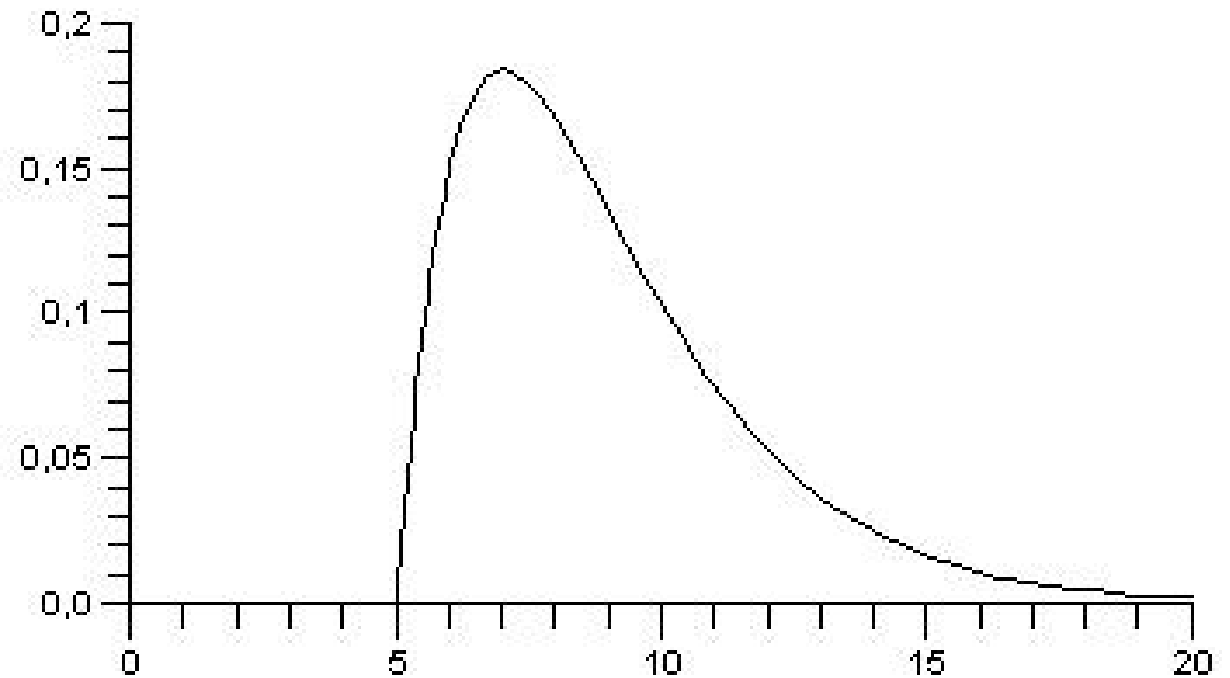} \vspace{-14pt}   
\\
 \small $ _\lambda$ ~~~~&  $ _\lambda$~~~~~
\end{tabular}
\caption{Left: example of a small-particle population mass distribution.  
Right: example of a large-particle population mass distribution.}
\label{Fig2} 
\end{figure}

In cases when sediment cohesion is insured by a biological factor with an impact
on the particles' mass, it is not reasonable to use a mass-preserving model. 
It is preferable to use a model able to reproduce the fact that for two sediment-particle populations
issued from the same initial sediment-particle population - one made up of small particles and the other 
of large particles - the total mass of the second population is 
greater than the total mass of the first. 
Figure \ref{Fig2} shows two mass distributions with respect to $\lambda$. 
Their total mass is not the same. 
A non-preserving mass model will be able to generate mass distributions of
those shapes, from the same initial mass distribution.    

\subsubsection{Operator building and properties}
To achieve the goal presented above, a decreasing function $\f(\lambda)$ is introduced. Then, using
\begin{gather}
\label{RM11} 
\G(\ffl(t,\pos),\rhc(t,\pos,.),\lambda) = -\frac{1}{\tr}
\bigg( \rhc(t,\pos,\lambda) 
-\frac{\ds \int_{\Lambda} \rhc(t,\pos,\lambda')\f(\lambda') \, d\lambda' }{\f(\lambda)}
\De(\ffl(t,\pos),\lambda)\bigg)
\end{gather}
in (\ref{GI4}) pushes the solution toward 
\begin{gather}
\frac{\ds \int_{\Lambda} \rhc(t,\pos,\lambda')\f(\lambda') \, d\lambda' }{\f(\lambda)}
\De(\ffl(t,\pos),\lambda),
\end{gather}
with a relaxation time $\tr$ and does not influence the evolution of 
\begin{gather}
\label{RM13} 
\int_{\Lambda} \rhc(t,\pos,\lambda)\f(\lambda) \, d\lambda .
\end{gather}
\subsubsection{Example of function $f$ for one-dimensional particles}
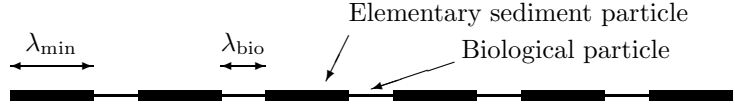
\begin{figure}
\setlength{\unitlength}{1mm}
\centering
\begin{picture}(95,20)(0,10)
\linethickness{1.2mm}
\multiput(0,10)(17,0){6}{\line(1,0){11}}
\linethickness{0.3mm}
\multiput(11,10)(17,0){5}{\line(1,0){6}}
\linethickness{0.1mm}
\put(0,14){\vector(1,0){11}}
\put(11,14){\vector(-1,0){11}}
\put(2,16){$\lm$}
\put(28,14){\vector(1,0){6}}
\put(34,14){\vector(-1,0){6}}
\put(28,16){$\lb$}
\put(45,20){Elementary sediment particle}
\put(45,18){\vector(-1,-2){3}}
\put(60,15){Biological particle}
\put(59,14.5){\vector(-3,-1){11}}
\end{picture}
\caption{Example of particle resulting from the aggregation of 
6 elementary sediment particles and 5 biological particles.}\label{figadisc}
\end{figure}

A suitable function $\f$ for one-dimensional particles aggregated by means of biological
factors may be built as follows. If the particles are the result of assemblies of 
elementary sediment particles with length $\lm$, joined together with biological
particles of length $\lb$, the length of a given particle is 
\begin{gather}
\label{RM14} 
\lambda = n \lm + (n-1) \lb,
\end{gather}
for a given $n$. An example of such a particle, with $n=6$, is represented in 
figure \ref{figadisc}.
Obviously, $n$ can be expressed in terms of $\lambda$, $\lm$ and $\lb$.
Indeed, since $(n-1) \lm + (n-1) \lb=\lambda-\lm$, the following formula are true:
\begin{gather}
\label{RM15} 
n-1 = \frac{\lambda-\lm}{\lm+\lb} \text{ ~ and ~ } n = \frac{\lambda+\lb}{\lm+\lb}.
\end{gather}
The following quantities may also be computed:
\begin{gather}
\label{RM16} 
\begin{aligned}
&\theta(\lambda) = \frac{n\lm}{\lambda} = \frac{(\lambda+\lb)\lm}{(\lm+\lb)\lambda}\leq 1,\\
& (1-\theta(\lambda)) = \frac{(n-1)\lb}{\lambda} = \frac{(\lambda-\lm)\lb}{(\lm+\lb)\lambda}\leq 1.
\end{aligned}
\end{gather}
Those quantities are the lineic proportions of elementary sediment particles and of
biological particles in a particle of length $\lambda$.
As a matter of fact, the mass $\mass(\lambda)$ of a particle of
length $\lambda$ may be expressed as
\begin{gather}
\label{RM17} 
\mass(\lambda) = \theta(\lambda) \lambda \mm + (1- \theta(\lambda)) \lambda \mb
= \frac{(\lambda+\lb)\lm}{(\lm+\lb)}\mm +\frac{(\lambda-\lm)\lb}{(\lm+\lb)}\mb,
\end{gather}
where $\mm$ is the lineic mass density of elementary sediment particles, and
$\mb$ is the lineic mass density of biological particles.

Beside this, when a particle of length $\lambda$ joins up with another of length
$\lambda'$, the result is a particle of length 
\begin{gather}
\label{RM18} 
\lambda'' =\lambda + \lambda' + \lb,
\end{gather}
which satisfies:
\begin{gather}
\label{RM19} 
\theta(\lambda'')\lambda'' =\theta(\lambda)\lambda + \theta(\lambda')\lambda',
\end{gather}
since 
\begin{gather}
\label{RM19.1} 
\frac{(\lambda''+\lb)\lm}{(\lm+\lb)} = \frac{(\lambda+\lb)\lm}{(\lm+\lb)} +\frac{(\lambda'+\lb)\lm}{(\lm+\lb)}.
\end{gather}
Multiplying (\ref{RM19}) by $\mm$ and rewriting it as 
\begin{multline}
\label{RM20} 
\frac{\theta(\lambda'')\lambda''\mm}{\theta(\lambda'') \lambda'' \mm + (1- \theta(\lambda'')) \lambda'' \mb}
\mass (\lambda'')=
\\
\frac{\theta(\lambda)\lambda\mm}{\theta(\lambda) \lambda \mm + (1- \theta(\lambda)) \lambda \mb}
\mass (\lambda)
+\frac{\theta(\lambda')\lambda'\mm}{\theta(\lambda') \lambda' \mm + (1- \theta(\lambda')) \lambda' \mb}
\mass (\lambda'),
\end{multline}
it may be deduced that defining
\begin{gather}
\label{RM22} 
\f(\lambda) = 
\frac{\theta(\lambda)\mm}{\theta(\lambda)  \mm + (1- \theta(\lambda)) \mb} =
\frac{1}{\ds 1+ \frac{1- \theta(\lambda)}{\theta(\lambda)} \frac{\mb}{\mm}} =
\frac{1}{\ds 1+ \frac{\lambda-\lm}{\lambda+\lb} \frac{\lb\mb}{\lm\mm} },
\end{gather}
equality 
\begin{gather}
\label{RM23} 
\f(\lambda'') m(\lambda'') = \f(\lambda) m(\lambda)  + \f(\lambda') m(\lambda') ,
\end{gather}
is true.

As a consequence, if a sediment-particle population is made, for each $n\in\nit$,
of $N_1(\lambda(n))$ particles of length $\lambda(n)= n \lm + (n-1)\lb$, and if 
a second population is generated by joining and breaking down particles of the first
population, leading  to $N_2(\lambda(n))$ particles of length $\lambda(n)$, for each $n\in\nit$,
then the following link between sequences $(N_1(\lambda(n)))_{n\in\nit}$ 
and $(N_2(\lambda(n)))_{n\in\nit}$
holds
\begin{gather}
\label{Suppl1} 
\sum_{n\in\nit}N_1(\lambda(n))\, f(\lambda(n))\, \mass(\lambda(n)) = 
\sum_{n\in\nit}N_2(\lambda(n))\, f(\lambda(n))\, \mass(\lambda(n)).
\end{gather}
Since the mass measure, on $\Lambda$, of the two sediment-particle populations are
\begin{gather}
R_1 = \sum_{n\in\nit}N_1(\lambda(n))\,  \mass(\lambda(n)) 
\; \delta_{\lambda=\lambda(n)},
\\
R_2 = \sum_{n\in\nit}N_2(\lambda(n))\,  \mass(\lambda(n)) 
\; \delta_{\lambda=\lambda(n)},
\label{Supp2l} 
\end{gather}
where $\delta_{\lambda=\lambda(n)}$ stands for the Dirac mass located in
$\lambda(n)$,
formula (\ref{Suppl1}) expressed also as
\begin{gather}
\label{Suppl3} 
<R_1,f> = <R_2,f> .
\end{gather}
Now, if measures $R_1$ and $R_2$ are replaced (or approached) by
$\rhc_1(\lambda) \, d\lambda$ and $\rhc_2(\lambda) \, d\lambda$, 
formula (\ref{Suppl3}) yields 
\begin{gather}
\label{Suppl4} 
\int_{\Lambda} \rhc_1(\lambda)f (\lambda)\, d\lambda = 
\int_{\Lambda} \rhc_2(\lambda)f(\lambda) \, d\lambda.
\end{gather}

Then property (\ref{RM13}) leads to the conclusion that choosing $\f$ 
defined by (\ref{RM22}) in the definition (\ref{RM11})
of $\G$ induces a behavior of $\rhc$ with respect to $\lambda$ in accordance with
sediment that aggregates because of biological factors as described above. 

~

In cases when $\mb<\mm$ and $\lb<<\lm$, since $\lambda > \lm$, it is also clear
that $\lb<<\lambda$. Then, $\f (\lambda) $ may be approached by:
\begin{gather}
\label{RM24} 
1- \frac{\lambda-\lm}{\lambda+\lb} \frac{\lb\mb}{\lm\mm} \sim
1- \frac{\lambda-\lm}{\lambda} \frac{\lb\mb}{\lm\mm}
\end{gather}

\subsubsection{Example of function $f$ for any-dimensional particles} For particles with dimension $d$ (which may not be an integer) a suitable function 
may be built using a generalization of the above considerations.
As previously, the considered particles are the result of assemblies of elementary 
sediment particles stuck together with biological particles. A particle with characteristic
length $\lambda$ is considered to be made of a proportion $\theta(\lambda)$ of
elementary sediment particles, and of a proportion $1-\theta(\lambda)$ of biological
particles. The mass $\mass(\lambda)$ of such particles is then
\begin{gather}
\label{RM25} 
\mass(\lambda) = \theta(\lambda) \lambda^d \mm + (1- \theta(\lambda)) \lambda^d \mb,
\end{gather}
where $\mm$ and $\mb$ are linked with a $d-$dimensional mass density of  elementary 
sediment particles and a $d-$dimensional mass density of biological particles.
When a particle of characteristic length $\lambda$ joins with another of
characteristic length $\lambda'$, it generates a particle of characteristic length
$\lambda''$ such that  
\begin{gather}
\label{RM26} 
\theta(\lambda''){\lambda''}^d =\theta(\lambda)\lambda^d + \theta(\lambda'){\lambda'}^d,
\end{gather}
or  
\begin{multline}
\label{RM27} 
\frac{\theta(\lambda''){\lambda''}^d\mm}
        {\theta(\lambda'') {\lambda''}^d \mm + (1- \theta(\lambda'')) {\lambda'' }^d\mb}
\mass (\lambda'')=
\\
\frac{\theta(\lambda)\lambda^d\mm}
        {\theta(\lambda) \lambda^d \mm + (1- \theta(\lambda)) \lambda^d \mb}
\mass (\lambda)
+\frac{\theta(\lambda'){\lambda'}^d\mm}
          {\theta(\lambda'){ \lambda' }^d\mm + (1- \theta(\lambda')){ \lambda'}^d \mb}
\mass (\lambda').
\end{multline}
Then choosing 
\begin{gather}
\label{RM28} 
\f(\lambda) = 
\frac{\theta(\lambda)\mm}{\theta(\lambda)  \mm + (1- \theta(\lambda)) \mb},
\end{gather}
insures the following
\begin{gather}
\label{RM29} 
\f(\lambda'') m(\lambda'') = \f(\lambda) m(\lambda)  + \f(\lambda') m(\lambda') ,
\end{gather}
leading to the conclusion that choosing $\f(\lambda)$ defined by (\ref{RM28})
induces  a behavior of $\rhc$ in accordance with the behavior of the distribution of
$d-$dimensional particles, which aggregate by means of a biological factor.

\section{Statistical aggregation and fragmentation models}
\label{SAF}
In this section, aggregation and fragmentation models usually used in chemical 
engineering and colloid sciences and often referenced as "Monte-Carlo Simulation"
(see 
 Gardner \& Theis \cite{GarTh},
Spilman \& Levenspiel \cite{SpiLe},
Daniels \& Hughes  \cite{DaHu},
Meakin  \cite{Meakin},
Liffman \cite{Liffman},
Shah \emph{et al.} \cite{ShaDoBo},
Das \cite{Das},
Spouge \cite{Spouge},
and Van Peborgh \&  Hounslow \cite{PeHo})
are adapted to the framework presented in this paper.
\subsection{Statistical aggregation and fragmentation considerations}
Summarizing the ideas used in the above-cited references, and revisiting them with 
a viewpoint inspired from the Boltzmann equation context (see Cercignani  \cite{cercignani})
we are led to the following reasoning.
The variable 
$\lambda\in \Lambda=[\lm,+\infty)$ stands for the particle size which is minimized by the 
size $\lm$ of elementary sediment particles, and $\mass(\lambda)$ for the mass of particles
of size $\lambda$. If the particles are $d-$dimensional,
\begin{gather}
\label{SAF1} 
\mass(\lambda) = \nd\, \lambda^d,
\end{gather}
for a constant $\nd$ depending on $d$.

Aggregation is described introducing a transition probability $\Ba(\ffl,\lambda,\lambda')$,
which depends on the fluid field $\ffl$. By definition, $\Ba(\ffl,\lambda,\lambda')$ is the probability
that two particles, one of size $\lambda$ and one of size $\lambda'$, being in the same place 
in fluid conditions $\ffl$, aggregate in a unit of time. 
It has of course the following property
\begin{gather}
\label{SAF2.1} 
\Ba(\ffl,\lambda,\lambda') = \Ba(\ffl,\lambda',\lambda),
\end{gather}
for any $\ffl$, $\lambda$ and $\lambda'$.
(Expressions of transition probability $\Ba$, based on physical principles, may be found
in Kim \emph{et al.} \cite{KimLeeJooLee}.)

Once aggregated, the two particles
give a particle of size $\lambda'' = \lnppsum$ which is such that 
\begin{gather}
\label{SAF2} 
\mass(\lambda'') = \mass(\lambda)+\mass(\lambda') .
\end{gather}

Fragmentation is described using $\Bf(\ffl,\lambda)$, which is the probability that a given
particle of size $\lambda$, in a place with fluid conditions $\ffl$, fragments in a unit 
of time.  Then, $\Be(\ffl,\lambda,\lambda')$ is the probability density function, with respect
to variable $\lambda'$, that a particle of size $\lambda$ which fragments gives a particle of size
$\lambda'\leq\big( \lambda^d/2 \big)^{1/d}$ (and another of size 
$\lambda'' =\lnpp \geq\big( \lambda^d/2 \big)^{1/d}$). By definition, $\Be$ has the following
properties
\begin{gather}
\label{SAF3} 
\begin{gathered}
 \Be(\ffl,\lambda,\lambda') = 0 \text{ if } \lambda' > \big( \lambda^d/2 \big)^{1/d},
 \\
 \int_{\lambda'\in \Lambda} \Be(\ffl,\lambda,\lambda') \, d\lambda' =
\int_{\lambda'\leq\big( \lambda^d/2 \big)^{1/d}} \Be(\ffl,\lambda,\lambda') \, d\lambda' = 1.
\end{gathered}
\end{gather}
Denoting by $\Bet(\ffl,\lambda,\lambda'')$ the probability density function, with respect
to variable $\lambda''$, that a particle of size $\lambda$ which fragments gives a particle of size
$\lambda''$ with $\big( \lambda^d/2 \big)^{1/d}\leq \lambda'' \leq \lambda$,
and by $\Ld(\lambda') = \lnpp$, $\Lr(\lambda'') = \lnss$, for any set 
$\omega\subset [\big( \lambda^d/2 \big)^{1/d}, \lambda]$, $\Be$ and $\Bet$ are linked
by
\begin{gather}
\label{SAF3.1} 
\int_{\omega} \Bet(\ffl,\lambda,\lambda'') \, d\lambda'' = 
\int_{\Lr(\omega)} \Be(\ffl,\lambda,\lambda') \, d\lambda' ,
\end{gather}
since every time that a particle of size $\lambda'$ is created, another particle 
of size $\lambda''=\Ld(\lambda')$ is also created. On the other hand, since 
the derivative of $L$ is 
\begin{gather}
\label{SAF3.2} 
\Ld'(\lambda') = -\lnppm {\lambda'}^{d-1},
\end{gather} 
making the change of variables $\lambda'' \mapsto \lambda' = \Lr(\lambda'')$, we get
\begin{gather}
\label{SAF3.3} 
\int_{\omega} \Bet(\ffl,\lambda,\lambda'') \, d\lambda'' = 
\int_{\Lr(\omega)} \Bet(\ffl,\lambda, \lnpp) \, \lnppm {\lambda'}^{d-1} \, d\lambda' .
\end{gather} 
Hence, $\Be$ and $\Bet$ are linked by
\begin{gather}
\label{SAF3.4} 
\Bet(\ffl,\lambda, \lnpp) \, \lnppm {\lambda'}^{d-1} = \Be(\ffl,\lambda,\lambda'),
\end{gather}
or, since 
$\ds   \lnppm {\lambda'}^{d-1} = \lnssminv  {\lambda''}^{1-d}$ when $\lambda'' =\Ld(\lambda')$,
by
\begin{gather}
\label{SAF3.5} 
\Be(\ffl,\lambda, \lnss) \, \lnssm {\lambda''}^{d-1} = \Bet(\ffl,\lambda,\lambda'').
\end{gather}
\subsection{Building operator $\G$ }
Now, building an integral operator $\G$ to be used in (\ref{GI4}), which takes those 
facts into account, consists in considering that the evolution of $\rhc(t,\pos,\lambda)$
in the neighborhood of a given value of $\lambda$, is the result of the following factors: 
a loss due to the aggregation of particles of size $\lambda$ with others,
another loss due the fragmentation of 
particles of size $\lambda$, a gain due to the aggregation of particles smaller than 
$\lambda$ and another gain due to the fragmentation of particles bigger than $\lambda$.

~

Quantifying the fragmentation-linked loss consists in noting that the density of particles
at a given point $\pos$ and in a given size $\lambda$ is nothing but 
$\rhc(t,\pos,\lambda)/\mass(\lambda)$, and in considering that, per unit of time, the 
number of particles of size $\lambda$ to fragment is in proportion with the number of
present particles. Hence, 
\begin{gather}
\label{SAF4} 
\frac{\rhc(t,\pos,\lambda)}{\mass(\lambda)}\; \Bf(\ffl,\lambda),
\end{gather} 
is the density, with respect to variables $\pos$ and $\lambda$, of particles of size $\lambda$ which
fragment per unit of time in $\pos$ and at $t$, when fluid conditions are $\ffl$.
Then, the density of mass loss related to fragmentation is:
\begin{gather}
\label{SAF5} 
\mass(\lambda)\; \frac{\rhc(t,\pos,\lambda)}{\mass(\lambda)}\; \Bf(\ffl,\lambda)=
\rhc(t,\pos,\lambda) \,\Bf(\ffl,\lambda).
\end{gather}
In order to quantify loss linked to aggregation, it must be noted that the probability
of a particle of size $\lambda$ to aggregate with a particle of size $\lambda'$, over a unit
of time, is in direct proportion to the number of particles of size $\lambda$ and 
the number of particles of size $\lambda'$. Consequently, the density, 
with respect to variables $\pos$ and $\lambda$, of particles of size 
$\lambda$ which aggregate is:
\begin{gather}
\label{SAF6} 
\int_{\lambda'\in \Lambda} 
\frac{\rhc(t,\pos,\lambda)}{\mass(\lambda)}
\frac{\rhc(t,\pos,\lambda')}{\mass(\lambda')}\; \Ba(\ffl,\lambda,\lambda') 
\; d\lambda',
\end{gather}
and the associated mass loss density, with respect to variables $\pos$ and $\lambda$, is
\begin{gather}
\label{SAF7} 
\int_{\lambda'\in \Lambda} 
{\rhc(t,\pos,\lambda)} \,
\frac{\rhc(t,\pos,\lambda')}{\mass(\lambda')}\; \Ba(\ffl,\lambda,\lambda') 
\; d\lambda'.
\end{gather}

~

The sum of the gain due to fragmentation is in two parts.
The first part is the result of fragmentations, the smallest resulting particles 
of which are of size $\lambda$, and the density of which, with respect to $\pos$
and $\lambda$ can be written as: 
\begin{gather}
\label{SAF8} 
\mass(\lambda)
\int_{\lambda' > \lambda} 
\frac{\rhc(t,\pos,\lambda')}{\mass(\lambda')}\; 
\Bf(\ffl,\lambda') \,\Be(\ffl,\lambda',\lambda) 
\; d\lambda'.
\end{gather}
In order to understand (\ref{SAF8}), it has to be noted that 
$\rhc(t,\pos,\lambda')/\mass(\lambda')$ is the density, with respect to variables $\pos$ and
$\lambda'$,  of particles of 
size $\lambda'$; $\Bf(\ffl,\lambda')$ is the probability of a particle of size
$\lambda'$ to fragment, within a unit of time; and, $\Be(\ffl,\lambda',\lambda)$ is the probability
density function, with respect to variable $\lambda$, of a particle of size $\lambda'$ 
to produce a particle
of size $\lambda$ as its smallest resulting particle. Hence in (\ref{SAF8}),
the integral is the density, with respect to variables $\pos$ and $\lambda$,  
of particles produced at size $\lambda$ as the
smallest particule resulting from fragmentation.
Multiplying this by $\mass(\lambda)$ gives the associated mass density.
Because of (\ref{SAF3}), (\ref{SAF8}) may be written as
\begin{gather}
\label{SAF9} 
\int_{\lambda' \in \Lambda} \mass(\lambda)
\frac{\rhc(t,\pos,\lambda')}{\mass(\lambda')}\; 
\Bf(\ffl,\lambda') \,\Be(\ffl,\lambda',\lambda) 
\; d\lambda'.
\end{gather}

The second part is the result of fragmentations whose largest resulting
particles are of size $\lambda$. The density, with respect to variables $\pos$ and 
$\lambda$,  associated with this second part reads:
\begin{gather}
\label{SAF10} 
\int_{\lambda' \in \Lambda} \mass(\lambda)
\frac{\rhc(t,\pos,\lambda')}{\mass(\lambda')}\; 
\Bf(\ffl,\lambda') \,\Bet(\ffl,\lambda',\lambda) 
\; d\lambda',
\end{gather}
or, because of (\ref{SAF3.5}),
\begin{gather}
\label{SAF11} 
\int_{\lambda' \in \Lambda} \mass(\lambda)
\frac{\rhc(t,\pos,\lambda')}{\mass(\lambda')}\; 
\Bf(\ffl,\lambda') \,\Be(\ffl,\lambda', \lpnp) \, \lpnpm {\lambda}^{d-1}
\; d\lambda'.
\end{gather}

As a consequence of (\ref{SAF9}) and (\ref{SAF11}), it may be concluded
that the mass-gain density, with respect to $\pos$ and $\lambda$,
due to fragmentation result is 
\begin{gather}
\label{SAF12} 
\int_{\lambda' \in \Lambda} \mass(\lambda)
\frac{\rhc(t,\pos,\lambda')}{\mass(\lambda')}\; 
\Bf(\ffl,\lambda') 
\bigg(\Be(\ffl,\lambda',\lambda) 
+  \Be(\ffl,\lambda', \lpnp) \, \lpnpm {\lambda}^{d-1}\bigg)
\; d\lambda'.
\end{gather}

The aggregation-related mass-gain density, with respect to $\pos$ and $\lambda$, reads
\begin{gather}
\label{SAF13.0} 
\int_{\lambda'^d\leq \lambda^d/2}
\mass(\lambda) \frac{\rhc(t,\pos,\lambda')}{\mass(\lambda')}
\frac{\rhc(t,\pos,\lnpp)}{\mass(\lnpp)} \lnppm \lambda^{d-1} 
\Ba(\ffl,\lambda', \lnpp)
\, d\lambda',
\end{gather} 
or using (\ref{SAF2.1}),
\begin{gather}
\label{SAF13} 
\frac12 \int_{\lambda'^d\leq \lambda^d}
\mass(\lambda) \frac{\rhc(t,\pos,\lambda')}{\mass(\lambda')}
\frac{\rhc(t,\pos,\lnpp)}{\mass(\lnpp)} \lnppm \lambda^{d-1}
\Ba(\ffl,\lambda', \lnpp)
\, d\lambda'.
\end{gather} 

~

The form of the operator $\G$ which may be deduced from (\ref{SAF5}), (\ref{SAF7}), (\ref{SAF12})
and (\ref{SAF13}) is 
\begin{multline}
\label{SAF14} 
\G (\ffl,\rhc(t,\pos,.),\lambda)=
- \rhc(t,\pos,\lambda) \,\Bf(\ffl,\lambda)
- \int_{\lambda'\in \Lambda} 
{\rhc(t,\pos,\lambda)} \,
\frac{\rhc(t,\pos,\lambda')}{\mass(\lambda')}\; \Ba(\ffl,\lambda,\lambda') 
\; d\lambda'
\\
+ \int_{\lambda' \in \Lambda} \mass(\lambda)
\frac{\rhc(t,\pos,\lambda')}{\mass(\lambda')}\; 
\Bf(\ffl,\lambda') 
\bigg(\Be(\ffl,\lambda',\lambda) 
+  \Be(\ffl,\lambda', \lpnp) \, \lpnpm {\lambda}^{d-1}\bigg)
\; d\lambda'
\\
+ \frac12 \int_{\lambda'\leq \lambda}
\mass(\lambda) \frac{\rhc(t,\pos,\lambda')}{\mass(\lambda')}
\frac{\rhc(t,\pos,\lnpp)}{\mass(\lnpp)} 
\Ba(\ffl,\lambda', \lnpp)\lnppm \lambda^{d-1} 
\, d\lambda'.
\end{multline}
\subsection{Properties of $\G$}
If an operator $\G$ of the kind defined by (\ref{SAF14}) is chosen in equation (\ref{GI4}), 
then it does not influence the evolution of the total mass of sediment. 
In other words, a function $\rhc(t,\lambda)$,
not depending on $\pos$, solution to equation (\ref{RM6}) with $\G$ given by (\ref{SAF14}),
satisfies property (\ref{RM7}) since 
\begin{gather}
\label{SAF15} 
\int_{\lambda\in\Lambda} \G (\ffl,\rhc(t,.),\lambda)\, d\lambda =0.
\end{gather}

This may be seen by computing, on the one hand, the integral with respect to $\lambda$ 
of the third term of (\ref{SAF14}). Using (\ref{SAF1}), it gives
\begin{multline}
\label{SAF16.0} 
\int_{\lambda \in \Lambda}\int_{\lambda' \in \Lambda} \mass(\lambda)
\frac{\rhc(t,\lambda')}{\mass(\lambda')}\; 
\Bf(\ffl,\lambda') 
\bigg(\Be(\ffl,\lambda',\lambda) 
\\ ~~~~~~~~~~~~~~~~~~~~~~~~
+  \Be(\ffl,\lambda', \lpnp) \, \lpnpm {\lambda}^{d-1}\bigg)
\; d\lambda' d\lambda=
\\
\int_{\lambda \in \Lambda} \int_{\lambda' \in \Lambda}
\rhc(t,\lambda') \; 
\Bf(\ffl,\lambda') 
\frac{\lambda^d}{\lambda'^d}\Be(\ffl,\lambda',\lambda) 
\; d\lambda'd\lambda  ~~~~~~~~~~~~
\\~~~~~~~~~~~~~~~~~~~~~~~~
+ \int_{\lambda \in \Lambda} \int_{\lambda' \in \Lambda} 
\rhc(t,\lambda') \Bf(\ffl,\lambda') 
 \Be(\ffl,\lambda', \lpnp) \, \lpnpm \frac{\lambda^d}{\lambda'^d}{\lambda}^{d-1}
\; d\lambda'd\lambda.
\end{multline}
Making the change of variables $(\lambda,\lambda')\mapsto(\tilde\lambda,\tilde\lambda')$,
with $\tilde\lambda=\lpnp$ and $\tilde\lambda'=\lambda'$, in the last integral, since
$\lambda^{d-1} d\lambda'd\lambda = {\tilde\lambda}^{d-1} d\tilde\lambda' d \tilde\lambda$,
it gives
\begin{multline}
\int_{\lambda \in \Lambda} \int_{\lambda' \in \Lambda}
\rhc(t,\lambda')\Bf(\ffl,\lambda')
 \Be(\ffl,\lambda', \lambda) \,   {\lambda}^{1-d}\, \frac{\lpnpd}{\lambda'^d} {\lambda}^{d-1}
\; d\lambda'd\lambda =
\\
\int_{\lambda \in \Lambda} \int_{\lambda' \in \Lambda}
\rhc(t,\lambda')\Bf(\ffl,\lambda')
 \Be(\ffl,\lambda', \lambda) \,  \frac{\lpnpd}{\lambda'^d}
\; d\lambda'd\lambda,
\end{multline}
and (\ref{SAF16.0}) yields
\begin{multline}
\label{SAF16} 
\int_{\lambda \in \Lambda}\int_{\lambda' \in \Lambda} \mass(\lambda)
\frac{\rhc(t,\lambda')}{\mass(\lambda')}\; 
\Bf(\ffl,\lambda') 
\bigg(\Be(\ffl,\lambda',\lambda) 
\\ ~~~~~~~~~~~~~~~~~~~~~~~~
+  \Be(\ffl,\lambda', \lpnp) \, \lpnpm {\lambda}^{d-1}\bigg)
\; d\lambda' d\lambda=
\\
\int_{\lambda' \in \Lambda}
\rhc(t,\lambda') \; 
\Bf(\ffl,\lambda') 
\bigg( \int_{\lambda \in \Lambda}
\Be(\ffl,\lambda',\lambda) \; d\lambda   \bigg) d\lambda' =
\int_{\lambda' \in \Lambda}
\rhc(t,\lambda') \; 
\Bf(\ffl,\lambda') \; d\lambda'.
\end{multline}
Once this computation is done, it is obvious that the integral with
respect to $\lambda$ of the third term of (\ref{SAF14}) is the opposite of the integral with
respect to $\lambda$ of the first term.

On the other hand, integrating the last term of (\ref{SAF14}) with respect to $\lambda$, gives
\begin{gather}
\label{SAF17} 
\frac12 \int_{\lambda \in \Lambda}\int_{\lambda'\leq \lambda}
\mass(\lambda) \frac{\rhc(t,\lambda')}{\mass(\lambda')}
\frac{\rhc(t,\lnpp)}{\mass(\lnpp)}  
\Ba(\ffl,\lambda', \lnpp) \lnppm \lambda^{d-1}
\, d\lambda'd\lambda .
\end{gather}
which, considering the following change of variables  
$(\lambda,\lambda')\mapsto(\tilde\lambda,\tilde\lambda')=(\lnpp,\lambda')$, 
or equivalently $(\tilde\lambda,\tilde\lambda') \mapsto(\lambda,\lambda')=
(({\tilde\lambda}^d + {{\tilde\lambda'}}~^{\hspace{-3pt}d})^{1/d},\tilde\lambda')$
for which
$\lambda^{d-1} d\lambda'd\lambda = {\tilde\lambda}^{d-1} d\tilde\lambda' d \tilde\lambda$,
yields
\begin{multline}
\label{SAF18} 
\frac12 \int_{\lambda \in \Lambda}\int_{\lambda' \in \Lambda}
\rhc(t, \lambda)\rhc(t, \lambda')
\frac{\mass(\lnppsum)}{\mass(\lambda)\mass(\lambda' )}\Ba(\ffl,\lambda', \lambda)
 \lambda^{1-d} \lambda^{d-1}
\, d\lambda'd\lambda =
\\
\frac12 \int_{\lambda \in \Lambda}\int_{\lambda' \in \Lambda}
\rhc(t, \lambda)\rhc(t, \lambda')
\Big(  \frac{1}{\mass(\lambda)} + \frac{1}{\mass(\lambda')} \Big) 
\Ba(\ffl,\lambda', \lambda)\, d\lambda'd\lambda =
\\
\int_{\lambda \in \Lambda}\int_{\lambda' \in \Lambda}
\frac{\rhc(t, \lambda)\rhc(t, \lambda')}{\mass(\lambda)}
\Ba(\ffl,\lambda', \lambda)\, d\lambda'd\lambda,
\end{multline}
where the first equality is obtained using (\ref{SAF1}), and the second using  (\ref{SAF2.1}).
As the last quantity in (\ref{SAF18}) is nothing but the opposite of the integral, with respect
to $\lambda$, of the second term of (\ref{SAF14}), (\ref{SAF15}) is true.

\section{On numerical applications}
\label{NUM}
In order to build numerical methods approximating equation (\ref{GI4}), it must be noted
that it is possible to make a splitting in time. For a given small time step $\Delta t$, this 
splitting routine consists, knowing an approximation
$\rhca(s, .,.)$ of $\rhc(s, .,.)$ at time $s$, 
in computing first $\rhcb(s+\Delta t, .,.)$, which is an approximation
of the solution $\rhct$ to
\begin{gather}
\begin{aligned}
\label{NUM1}  
&\fracp{\rhct}{t} 
+ U(\ffl,\lambda)\fracp{\rhct}{x}
+ V(\ffl,\lambda)\fracp{\rhct}{y}
+ \big(W(\ffl,\lambda)-\ssvel(\lambda,\rrr)\big)\fracp{\rhct}{z}
\\
& ~~~~~~~~~~~~~~~~~~~~~~~~~~
-\bigg(
\fracp{\Big(\mu(\ffl,\lambda)\ds\fracp{\rhct}{x}\Big)}{x}
+ \fracp{\Big(\mu(\ffl,\lambda)\ds\fracp{\rhct}{y}\Big)}{y}
+ \fracp{\Big(\nu(\ffl,\lambda)\ds\fracp{\rhct}{z}\Big)}{z}
\bigg)
= 0,
\\
&\rhct(s, .,.)=\rhca(s, .,.),
\end{aligned}
\end{gather}
at time $s+\Delta t$, and then $\rhca(s+\Delta t, .,.)$ as an approximation of the solution $\rhch$ to
\begin{gather}
\begin{aligned}
\label{NUM2}  
&\fracp{\rhch}{t} = \G (\ffl,\rhch,\lambda),
\\
&\rhch(s, .,.)=\rhcb(s+\Delta t, .,.),
\end{aligned}
\end{gather} 
at time $s+\Delta t$.

After discretizing $\Lambda$ into $I$ subsets $\Lambda_i$, an approximation of 
(\ref{NUM1}) is made of a collection of $I$ advection-diffusion equations (one for
each $\Lambda_i$), which are not mutually dependent.  An approximation of each
of these advection-diffusion equations may be computed using any usual 
numerical advection-diffusion solver. 
Concerning the computation of an approximation of the solution to (\ref{NUM2}),
which is the focus here, once the fluid field $\ffl$ is known and the position in space 
is discretized, it comes down to computing a collection of approximations of
$\rhc(s+\Delta t, \lambda)$, solution to 
\begin{gather}
\begin{aligned}
\label{NUM3}  
&\fracp{\rhc}{t} = \G (\rhc,\lambda),
&\rhc(s, \lambda)=\rhc_0(\lambda),
\end{aligned}
\end{gather} 
(forgetting the dependence in $\ffl$ and $\pos$) for given functions $\rhc_0$.

~

In the case when $\G$ is given by (\ref{SAF14}),
one way to proceed would be to follow a Monte-Carlo method 
(see Lapeyre, Pardoux \& Sentis \cite{LaPaSen}) to approximate
the integrals involved within the definition of $\G$.
Proceeding in this direction would also lead to numerical methods of the types
used in 
Spilman \& Levenspiel \cite{SpiLe},
Daniels Hughes  \cite{DaHu},
Meakin  \cite{Meakin},
Liffman \cite{Liffman},
Shah \emph{et al.} \cite{ShaDoBo},
Das \cite{Das},
Spouge \cite{Spouge},
Van Peborgh \&  Hounslow \cite{PeHo}
and  Kim \emph{et al.} \cite{KimLeeJooLee}
 
There is another way in which an approximated solution of (\ref{NUM3}) with $\G$
given by (\ref{SAF14}) can be built. This way consists in building a discrete
operator $\Gov$ from $\G$, without breaking the structure yielding property (\ref{SAF15}).
As a matter of fact, a discrete version of property (\ref{SAF15}) may be written for
the discrete operator $\Gov$. The method followed does have a relation to the method
set out in 
Buet  \cite{buet},
Rogier \& Schneider  \cite{rogier/schneider},
Degond \& Lucquin \cite{degond/lucquin:1992} and
Fr\'enod \& Lucquin  \cite{frenod/lucquin}
in the contexts of Boltzmann and Fokker-Plank equations.
\subsection{Discrete operator building }
First, as mentioned previously, $\Lambda$ is discretized into $I$ subsets 
$\Lambda_i=[\lambda_i, \lambda_{i+1})$ such that $\Lambda_i\cap\Lambda_j =\emptyset$
and $\Lambda =\cup_{i= 1}^{I}\Lambda_i.$ Then, denoting by $|\Lambda_i|$ the measure of
$\Lambda_i$, the following numbers are defined:
\begin{align}
&\Bfli = \frac{1}{|\Lambda_i|} \int_{\ldt\in\Lambda_i} \Bf(\ldt) \, d\ldt, 
\text{ for } i=1,\dots, I,
\label{NUM4} 
\\
& \Blai = \frac{1}{|\Lambda_i||\Lambda_j|} \int_{\ldt\in\Lambda_i}\int_{\ldtp\in\Lambda_j}
\frac{1}{\mass(\ldtp)} \Ba(\ldt,\ldtp)  \, d\ldtp d\ldt,
\text{ for } i \text{ and } j =1,\dots,  I,
\label{NUM5} 
\\ \nonumber
&\Bgfi =  \frac{1}{|\Lambda_i||\Lambda_j|}   \int_{\ldt\in\Lambda_i}\int_{\ldtp\in\Lambda_j}
\Bf(\ldtp)\frac{\mass(\ldt)}{\mass(\ldtp)} 
\bigg(\Be(\ldtp,\ldt) +
\\&~~
\Be(\ldtp, \lptnpt) \lptnptm \ldt^{d-1}  \lptnptm
\bigg)  \, d\ldtp d\ldt,
\text{ for } i \text{ and } j = 1,\dots,  I,
\label{NUM6} 
\\ \nonumber
&\Bgai =  \frac{1}{|\Lambda_i|}   \int_{\ldt\in\Lambda_i}   
\frac{1}{|\Lambda_{jk}(\lambda)|} 
\int_{\ldtp\in\Lambda_{jk}(\lambda)
 }
\frac{\mass(\ldt)}{\mass(\ldtp)\mass(\lnptpt)}
\\ &~  \hspace{2.3cm}
\Ba(\ldt, \lnptpt) \lnptptm \, d\ldtp d\ldt , \text{ for } i = 1,\dots, I, ~ k \leq i \text{ and } l  \leq i,
\label{NUM7} 
\end{align}
where, in equation (\ref{NUM7}), $\Lambda_{jk}(\lambda) =$
$(\big(\lambda^d-\lambda_{l+1}^d\big)^{1/d},\big(\lambda^d-\lambda_{l}^d\big)^{1/d}] \cap \Lambda_k$, which reads also $\Lambda_{jk}(\lambda) =$
$\{  \ldtp\in\Lambda_k, \lnptpt\in\Lambda_l\}$.
The following functions are also defined:
\begin{gather}
\label{NUM8} 
\begin{aligned}
&\Bfl (\lambda) = \sum_{i=1}^{I} \Bfli  ~\mathds{1}_{ \hspace{-1pt}\Lambda_i}\hspace{-1pt}(\lambda),
\\
& \Bla(\lambda,\lambda') = \sum_{i=1}^{I} \sum_{j=1}^{I} \Blai
   ~\mathds{1}_{ \hspace{-1pt}\Lambda_i}\hspace{-1pt}(\lambda)
   ~\mathds{1}_{ \hspace{-1pt}\Lambda_j}\hspace{-1pt}(\lambda'),
\\
& \Bgf(\lambda,\lambda') = \sum_{i=1}^{I}  \sum_{j=1}^{I} \Bgfi
   ~\mathds{1}_{ \hspace{-1pt}\Lambda_i}\hspace{-1pt}(\lambda)
   ~\mathds{1}_{ \hspace{-1pt}\Lambda_j}\hspace{-1pt}(\lambda'),
\\
& \Bga(\lambda,\lambda',\lambda^*) = \sum_{i=1}^{I}  \sum_{k=1}^{i}  \sum_{l=1}^{i} \Bgai
   ~\mathds{1}_{ \hspace{-1pt}\Lambda_i}\hspace{-1pt}(\lambda)
   ~\mathds{1}_{ \hspace{-1pt}\Lambda_j}\hspace{-1pt}(\lambda')
  ~\mathds{1}_{ ((\lambda^d-\lambda_{l+1}^d)^{1/d},(\lambda^d-\lambda_{l}^d)^{1/d}] }
       \hspace{-1pt}(\lambda^*),
\end{aligned}
\end{gather}
where $\mathds{1}_{ \hspace{-1pt}\Lambda_i}$ stands for the characteristic function
of set $\Lambda_i$.

With those definitions at hand, operator $\Gov(\rhc(t,.),\lambda)$ is set as
\begin{multline}
\label{NUM9} 
\Gov(\rhc(t,.),\lambda) =
- \rhc(t,\lambda) \,\Bfl(\lambda)
- \int_{\lambda'\in \Lambda} 
 \rhc(t,\lambda)  \rhc(t,\lambda')  \; \Bla(\ffl,\lambda,\lambda') \; d\lambda' 
 \\
 +  \int_{\lambda'\in \Lambda} 
 \rhc(t,\lambda') \Bgf(\lambda,\lambda')  \; d\lambda'
 + \int_{\lambda'\in \Lambda}  \int_{\lambda* \in \Lambda}
 \rhc(t,\lambda')  \rhc(t,\lambda^*)  \Bga(\lambda,\lambda',\lambda^*)  \; d\lambda' d\lambda^*.
 \end{multline}
  If $\rhc(t,.)$ is constant on every $\Lambda_i$, 
 with worth $\rhc^i(t)$ or, in other words, if
 \begin{gather}
\label{NUM10.0} 
\rhc(t,\lambda) =  \sum_{i=1}^{I}  \rhc^i(t)
  ~\mathds{1}_{ \hspace{-1pt}\Lambda_i}\hspace{-1pt}(\lambda),
\end{gather}
then $\Gov(\rhc(t,.),\lambda)$ has the following expression
\begin{multline}
\label{NUM10} 
\Gov(\rhc(t,.),\lambda) =
- \sum_{i=1}^{I} \Bfli  \rhc^i(t) ~\mathds{1}_{ \hspace{-1pt}\Lambda_i}\hspace{-1pt}(\lambda)
- \sum_{i=1}^{I} \sum_{j=1}^{I} \Blai  \rhc^i(t)  \rhc^j(t) 
   ~\mathds{1}_{ \hspace{-1pt}\Lambda_i}\hspace{-1pt}(\lambda)
 \\
 +  \sum_{i=1}^{I} \sum_{j=1}^{I} \Bgfi  \rhc^j(t)
  ~\mathds{1}_{ \hspace{-1pt}\Lambda_i}\hspace{-1pt}(\lambda)
 + \sum_{i=1}^{I}  \sum_{k=1}^{i}  \sum_{l=1}^{i} \Bgai   \rhc^k(t)  \rhc^l(t)
 ~\mathds{1}_{ \hspace{-1pt}\Lambda_i}\hspace{-1pt}(\lambda) .
 \end{multline}

 \subsection{Discrete operator properties}
 By its construction, operator $\Gov(\rhc(t,.),\lambda)$ defined by (\ref{NUM9}) is close to 
 operator $\G(\rhc(t,.),\lambda)$ defined by (\ref{SAF14}), as soon as  $\rhc$ is regular enough.
  
 On the other hand, as a direct consequence of (\ref{NUM10}), if $\rhc(t,.)$ is constant on every 
 $\Lambda_i$,   then $\Gov(\rhc(t,.),\lambda)$  
 is constant on every $\Lambda_i$ and it is easy to see, as a consequence
 of the building of $\Gov$, that
 \begin{gather}
\int_{\Lambda_i} \Gov (\rhc(t,.),\lambda) \, d\lambda =
\int_{\Lambda_i} \G (\rhc(t,.),\lambda) \, d\lambda ,
\end{gather}
for $i=1,\dots,I$.  Then from equality (\ref{SAF15}), it may be deduced that 
\begin{gather}
\label{NUM11} 
\int_{\lambda\in\Lambda} \Gov (\rhc(t,.),\lambda)\, d\lambda =0.
\end{gather}

~

Hence it may be concluded that, if $\breve\rhc_0$ is constant on every $\Lambda_i$, the
solution to 
\begin{gather}
\begin{aligned}
\label{NUM3333}  
&\fracp{\breve\rhc}{t} = \Gov (\breve\rhc,\lambda),
& \breve\rhc(s, \lambda)= \breve\rhc_0(\lambda),
\end{aligned}
\end{gather} 
is also constant on every $\Lambda_i$ and satisfies
\begin{gather}
\label{NUM13} 
\int_{\lambda\in\Lambda} \breve\rhc(t,\lambda)\, d\lambda
\text{ is constant along time. }
\end{gather}

Consequently, a good way to build a mass-preserving numerical scheme to approximate
 (\ref{NUM3}) consists in approximating $\rhc_0$ by $\breve\rhc_0$ defined by
 \begin{gather}
\label{NUM14.00} 
\breve\rhc_0 (\lambda) = \sum_{i=1}^{I}\rhc_0^i \,
\mathds{1}_{ \hspace{-1pt}\Lambda_i}\hspace{-1pt}(\lambda)
\text{ ~ with ~ } \rhc_0^i = \frac{1}{|\Lambda_i|} \ \int_{\ldt\in\Lambda_i}
\rhc_0(\ldt) \, d \ldt,
\end{gather}
and in approximating $\rhc(s+\Delta t,\lambda)$, the solution to (\ref{NUM3}) at time $s+\Delta t$,  
by $\breve\rhc(s+\Delta t, \lambda)$ defined as :
\begin{gather}
\label{NUM14} 
\breve\rhc(s+\Delta t, \lambda)  =  \breve\rhc_0 (\lambda) + \Delta t  \, \Gov (\breve\rhc_0 (.),\lambda).
\end{gather}
The result $\breve\rhc(s+\Delta t, \lambda)$ will be close to 
$\rhc(s+\Delta t,\lambda)$, constant on every $\Lambda_i$  
and will satisfy 
\begin{gather}
\label{NUM14.1} 
\int_{\lambda\in\Lambda} \breve\rhc(s+\Delta t,\lambda) \, d\lambda =
\int_{\lambda\in\Lambda} \breve\rhc_0 (\lambda) \, d\lambda= 
\sum_{i=1}^{I} |\Lambda_i| \, \rhc_0^i  = 
\int_{\lambda\in\Lambda} \rhc_0 (\lambda) \, d\lambda.
\end{gather}
\section{Perspectives}
This paper puts forward a framework designed to process the characteristic evolution 
of sediment particles being transported by the water column.
It gives simple examples of instantiations of this framework.

~

It opens the way for many interesting questions.

Concerning modeling, it would be useful to incorporate other sediment particle characteristics 
into the model, such as charge and fractal dimension, which seem important and are attentively
studied by colloid scientists. To do this, new spaces $\Lambda$ of larger dimension
have to be built, taking into account the input from colloid sciences.
 
Concerning mathematics, the existence of results for equations of the kind
\begin{gather}
\fracp{\rhc}{t} = \G (\ffl,\rhc,\lambda),
\end{gather} 
with $\G$ given by (\ref{RM3}), 
(\ref{RM11}) or (\ref{SAF14}) are an interesting challenge.

Finally, concerning numerics, the path explored in section \ref{NUM} has to be explored in greater depth,
and software has to be designed to test the accuracy of such schemes.

~

{\bf Acknowledgments - }The author thanks Juliette Bouchery for proofreading the manuscript.

\bibliographystyle{plain}
\bibliography{biblio}

\end{document}